\documentclass[12 pt]{amsart}

\newtheorem{theorem}{Theorem}[section]
\newtheorem{lemma}[theorem]{Lemma}

\newtheorem{corollary}[theorem]{Corollary}
\newtheorem{problem}[theorem]{Problem}

\theoremstyle{definition}

\newtheorem{example}[theorem]{Example}

\theoremstyle{remark}

\usepackage{amscd,amssymb}

\begin{document}
\baselineskip 15.5pt

\title{Coordinates, retracts and automorphisms}
\author{Yun-Chang Li}
                                         
\address{Department of Mathematics\\ The University of Hong Kong\\ Hong Kong SAR, China}
\email{liyunch@hku.hk\\ liyunch88@gmail.com}

\thanks
{The research of Yun-Chang Li  was partially supported by a
postgraduate studentship in the University of Hong Kong.}



\begin{abstract}

Let $K$ be a field of characteristic zero, $K[x,y]$ be the polynomial ring in two variables. Let $\phi=(f, g)$ be an endomorphism of $K[x,y]$. It is proved that if $\phi$ maps each coordinate to a generator of some proper retract, then it is an automorphism. As a corollary, the retract preserving problem is solved for both polynomial ring over $K$ and free algebra over an arbitrary field when $n=2$.

\bigskip

\end{abstract}

\maketitle

\section{Introduction}

\noindent Let $K$ be a field of characteristic zero, $P_n=K[x_1,\cdots, x_n]$ be the polynomial ring in $n$ variables and $\sigma$ be an endomorphism. By $\sigma$ "preserving" things of type $A$, we mean that $\sigma$ maps things of type $A$ to things of the same type. For some preserving problems of polynomial and free algebras, see \cite{100,500,700,800,900}. Recall that a subring $R$ of $P_n$ is a retract if there exists some idempotent endomorphism $\pi$ ($\pi$ is idempotent means $\pi^2=\pi$) such that $\pi(P_n)=R$. The endomorphism is called the corresponding retraction of $R$. For more information, please refer to \cite{100,200,300}. Then the corresponding "retract preserving problem" is raised naturally:

\

\begin{problem}[Jie-Tai, Yu]\label{p100}

\noindent Is an endomorphism that maps each proper retract to a proper  retract an automorphism ?

\end{problem}

\

\noindent For free algebras, the parallel problem is not true with the following simple counterexample.

\

\begin{example}

\noindent Let $F$ be an arbitrary field, $F\langle x,y\rangle$ be the free algebra in two variables and $\phi=(x, y+xy-yx)$ be an endomorphism. Then it is not an automorphism obviously since the leading terms are not algebraically dependent. According to \cite{200, 300}, if $R$ is a proper retract of $F\langle x, y\rangle$, then $R=F[r]$ for some $r\in F\langle x,y\rangle$. Let $\pi=(s(r),t(r))$ be the corresponding retraction, and then 
$$\pi(r)=r(s(r),t(r))=r.$$

\

\noindent Define $\phi(r)=r'$, and then $r'(x,y)=r(x,y)+c(x,y)$ where $c(x,y)$ is a commutator. Since $r(s(r), t(r))=r$, $r(s(r'), t(r'))=r'$. Hence if we define $\pi'=(s(r'), t(r'))$, then $\pi'(r')=r(s(r'), t(r'))+c(s(r'),t(r'))$. However, since $s(r')$ and $t(r')$ are algebraically dependent, $c(s(r'), t(r'))=0$, and hence $\pi'(r')=r'$ which implies that $F[r']$ is a retract with the corresponding retraction $\pi'$.

\end{example}

\

\noindent In this paper, a positive solution is given to Problem \ref{p100} when $n=2$. Following is the main theorem .

\

\begin{theorem}[Main]\label{t100}

\noindent Let $K$ be a field of characteristic zero, $K[x,y]$ be the polynomial ring in two variables and $\phi$ be an endomorphism that maps each coordinate to a generator of some proper retract. Then $\phi$ is an automorphism.

\end{theorem}

\

\begin{corollary}

\noindent Let $K$ be a field of characteristic zero. Then an endomorphism $\phi$ which preserves the proper retracts of $P_2$ is an automorphism.

\end{corollary}

\

\begin{proof}

\noindent Since each coordinate is also a generator of some proper retract, $\phi$ maps each coordinate to a generator of some proper retract, and hence it is an automorphism by Theorem \ref{t100}.

\end{proof}

\

\section{Proof of the Main Theorem}

\noindent Throughout this section, we always assume that $K$ is a field of characteristic zero. Let $\phi=(f,g)$ be an endomorphism which maps each coordinate to a generator of some proper retract. 

\

\begin{lemma}\label{l100}

\noindent Let $p$ be a coordinate and $p'=\phi(p)$. Let $\pi=(s(p'), t(p'))$ be some corresponding retraction of $K[p']$ and $\pi'=(s(z), t(z))$ where $s(z), t(z)\in K[z]$. Then $K[\pi'(f), \pi'(g)]=K[z]$. 

\end{lemma}

\

\begin{proof}

\noindent Since $\pi$ is the retraction, then $\pi(p')=p'$, or $p'(s(p'), t(p'))=p'$. Since 
$$p'(s(p'), t(p'))=p'(s(z), t(z))\mid_{z=p'},$$
$\pi'(p')$ has to be equal to $z$ since $p'\not\in K$. Hence $\pi'\circ\phi(p)=z$ which implies $K[\pi'(f), \pi'(g)]=K[z]$.

\end{proof}

\

\noindent By \cite{300}, since $\phi(x)=f$ is a generator of some proper retract, then there exists some automorphism $\sigma$ such that $\sigma(f)=x+y\cdot h_1(x,y)$. Let $\sigma(g)=y\cdot h_2(x,y)+h(x)$ and $\sigma'=(x, y-h(x))$, then $\sigma\circ\phi\circ\sigma'=(x+y\cdot h_1(x,y), y\cdot h_2(x,y))$. Obviously that $\sigma\circ\phi\circ\sigma'$ also maps each coordinate to a  generator of a proper retract and $\phi$ is an automorphism if and only if $\sigma\circ\phi\circ\sigma'$ is an automorphism. Hence we assume that $\phi=(x+y\cdot h_1(x,y), y\cdot h_2(x,y))$ for some $h_1, h_2\in K[x,y]$ and $h_2\not=0$. 

\

\begin{lemma}\label{l200}

\noindent To any $n$, there exists some $\pi_n=(s_n(z), t_n(z))$ where $\deg(s_n)\cdot\deg(t_n)>0$ and $\deg(s_n)>n\cdot\deg(t_n)$ such that 
$$K[\pi_n(f), \pi_n(g)]=K[z].$$

\end{lemma}

\

\begin{proof}

\noindent Assume not, and then there exists some positive integer $N$ such that if $\deg(s(z))\cdot\deg(t(z))>0$ and $K[f(s,t), g(s,t)]=K[z]$, then $\deg(s)\leq N\cdot\deg(t)$. 

\

\noindent Now consider the coordinate $y+(x+y^M)^2$ where $M>max\{ \deg(h_1)+2, N, 1+(N+1)\deg(h_1)\}$. By Lemma \ref{l100} there exists some $\pi=(s(z), t(z))$ such that $\pi\circ\phi(y+(x+y^M)^2)=z$ , or
$$\pi(yh_2+(x+yh_1+y^Mh_2^M)^2)=z.$$
If $t(z)$ is a constant, then $\pi(h_2)$ can not be a constant since if so, then $\pi(yh_2)$ is a constant, and hence $z-\pi(yh_2)=(\pi(x+yh_1+y^Mh_2^M))^2$ which is impossible. Hence $\deg(\pi(yh_2))\geq 1$ and $\deg(\pi(y^Mh_2^M))\geq M$. Since $K[s(z), t(z)]=K[z]$, then $s(z)=az+b$ where $a, b\in K$ and $a\not=0$ ($t(z)$ is a constant). Hence $\deg(\pi(yh_1))\leq \deg_x(h_1)<M$, and then 
$$\deg(\pi(yh_2+(x+yh_1+y^Mh_2^M)^2))=\deg(\pi(y^{2M}h_2^{2M}))>2,$$ 
which contradicts to $\pi(yh_2+(x+yh_1+y^Mh_2^M)^2)=z$. Hence $t(z)\not\in K$.

\

\noindent If $s(z)\in K$, then $t(z)=az+b$ where $a, b\in K$ and $a\not=0$. If $\pi(h_2)=0$, then 
$$\pi(yh_2+(x+yh_1+y^Mh_2^M)^2)=(\pi(x+yh_1))^2\not=z.$$ 
Hence $\pi(h_2)\not=0$, and then $\deg(\pi(yh_2))\geq 1$. Since $\deg(\pi(yh_1))\leq 1+\deg_y(h_1)<M$, 
$$\deg(\pi(x+yh_1+y^Mh_2^M))=\deg(\pi(y^Mh_2^M))\geq M.$$
Moreover, since $M\geq 2$, $\deg(\pi(y^{2M}h_2^{2M}))>\deg(\pi(yh_2))$, and hence 
$$\deg(\pi(yh_2+(x+yh_1+y^Mh_2^M)^2))=\deg(\pi(y^{2M}h_2^{2M}))>1$$
which contradicts to $\pi(yh_2+(x+yh_1+y^Mh_2^M)^2)=z$. Hence $s(z)\not\in K$.

\

\noindent Now assume $\deg(s)\cdot\deg(t)>0$. Also, $\pi(h_2)\not=0$ since $z$ is not a square. Then $\deg(\pi(yh_2))\geq 1$ and hence $\deg(\pi(y^Mh_2^M))\geq M\cdot\deg(t(z))$. 

\

\noindent By the condition, $\deg(s)\leq N\cdot\deg(t)$, and hence $\deg(\pi(x))<M\cdot \deg(t)$. Similarly,
$$\begin{array}{rcl}\deg(\pi(yh_1))&\leq&\deg(t)\cdot(1+\deg_y(h_1))+\deg(s)\cdot\deg_x(h_1)\\&\leq&\deg(t)+(\deg(s)+\deg(t))\deg(h_1)\\&\leq&\deg(t)\cdot(1+(N+1)\deg(h_1))\\&<&M\cdot\deg(t),\end{array}$$
and hence $\deg(\pi(yh_2+(x+yh_1+y^Mh_2^M)^2))=(\deg(\pi(x+yh_1+y^Mh_2^M))^2>1$ which contradicts. Hence to any positive integer $n$, there exists some $\pi_n=(s_n(z), t_n(z))$ where $\deg(s_n)\cdot\deg(t_n)>0$ and $\deg(s_n)>n\cdot\deg(t_n)$ such that 
$$K[\pi'(f), \pi'(g)]=K[z].$$

\end{proof}

\

Now we establish the lexicographic order on all monomials of $K[x,y]$ by $x>>y$ and denote the leading monomial of $l(x,y)$ by $v(l)$ for each polynomial $l(x,y)\in K[x,y]$.

\begin{lemma}\label{l300}

\noindent Let $\psi=(f', g')$ be an endomorphism of $K[x,y]$ that maps each coordinate to a generator of some proper retract. If there exists a sequence $(\pi_n=(s_n(z), t_n(z)))$ with $\deg(s_n)\cdot\deg(t_n)>0$ and $K[\pi_n(f'), \pi_n(g')]=K[z]$ such that $\deg(s_n)>n\cdot \deg(t_n)$, then one of $f', g'$ is of the form $a'y+b'$ where $a',b'\in K$, $a'\not=0$, or $v(f')$ and $v(g')$ are algebraically dependent with the one of a greater degree being a power of the other one.

\end{lemma}

\

\begin{proof}

\noindent If neither $f'$ nor $g'$ is of the form $a'y+b'$, then $x$ appears in both $f'$ and $g'$ since a polynomial of outer rank one can not be a generator of a proper retract if the degree is greater than 1. Assume $v(f')=y^ax^b$ and $v(g')=y^cx^d$ where $b\cdot d\not=0$. Then to any $y^ix^j\in supp(f')$, if assume $m_n=\deg(s_n)/\deg(t_n)$, we have
$$\frac{\deg(\pi_n(y^ix^j))}{\deg(\pi_n(y^ax^b))}=\frac{i+jm_n}{a+bm_n}$$
where $m_n>n$. Since $y^ax^b$ is the leading term, $j<b$ or $j=b, i<a$, and hence there exists some $N_{ij}$ such that when $n>N_{ij}$, $\deg(\pi_n(y^ix^j))<\deg(\pi_n(y^ax^b))$. Since there exists finitely many monomials in $supp(f')$, there exists some $N_1$ such that to any $n>N_1$ we have $\deg(\pi_n(f'))=\deg(\pi_n(y^ax^b))$. Similar to $f'$, there exists some $N_2$ such that to any $n>N_2$ we have $\deg(\pi_n(g'))=\deg(\pi_n(y^cx^d))$. Define $N_0=max\{N_1, N_2\}$, and then to any $n>N_0$ we have $\deg(\pi_n(f'))=\deg(\pi_n(y^ax^b))$ and $\deg(\pi_n(g'))=\deg(\pi_n(y^cx^d))$.

\

\noindent Since $K[\pi_n(f'), \pi_n(g')]=K[z]$, by the famous Abhyankar-Moh Theorem (\cite{400}, Main Theorem), the greater one of $\deg(\pi_n(f'))$ and $\deg(\pi_n(g'))$ is a multiple of the other one, and hence to any $n>N_0$, we also have the greater one of $\deg(\pi_n(f'))$ and $\deg(\pi_n(g'))$ is a multiple of the lower one. Without loss of generality, we can assume $b\geq d$, and then 
$$\lim_{n\rightarrow\infty}\frac{\deg(\pi_n(f'))}{\deg(\pi_n(g'))}=\frac{b}{d},$$
and if $b$ is not a multiple of $d$, then there exists some $n'$ such that $\deg(\pi_n(f'))$ is not a multiple of $\deg(\pi_n(g'))$ which contradicts. Hence $d\mid b$, and then 
$$\frac{\deg(\pi_n(f'))}{\deg(\pi_n(g'))}=\frac{a\deg(t_n)+b\deg(s_n)}{c\deg(t_n)+d\deg(s_n)}=k+\frac{a-ck}{c+dm_n}$$
where $k=b/d$. If $a-ck\not=0$, then $m_n$ can be chosen great enough such that $0<(a-ck)/(c+dm_n)<1$, and hence $\deg(\pi_n(f'))$ is no longer a multiple of $\deg(\pi_n(g'))$ which contradicts. Hence $a-ck=0$ which implies that $y^ax^b=(y^cx^d)^k$. 

\end{proof}

\

\begin{lemma}\label{l400}

\noindent Let $\psi=(f', g')$ be an endomorphism which maps each coordinate to a generator of some proper retract. If there exists some automorphism $\alpha$ such that $\psi\circ\alpha$ is of the form $(x+yh'_1, yh'_2)$, then there exists a sequence $(\pi_n=(s_n(z), t_n(z)))$ with $\deg(s_n)\cdot\deg(t_n)>0$ and $K[\pi_n(f'), \pi_n(g')]=K[z]$ such that $\deg(s_n)>n\cdot\deg(t_n)$.

\end{lemma}

\

\begin{proof}

\noindent By Lemma \ref{l200}, this sequence does exist for $\psi\circ\alpha$, and hence it is also a sequence satisfied to $\psi$. 

\end{proof}

\

\noindent \textbf{Proof for Main Theorem.}  Let $\phi=(f,g)$ be an endomorphism which maps each coordinate to a generator of some proper retract. Then it suffices to prove the special case of $f=x+yh_1$ and $g=yh_2$. By Lemma \ref{l200}, there exists a sequence $(\pi_n=(s_n(z), t_n(z)))$ with $\deg(s_n)\cdot\deg(t_n)>0$ and $K[\pi_n(f),\pi_n(g)]=K[z]$ such that $\deg(s_n)>n\cdot\deg(t_n)$, and by Lemma \ref{l300} there exists some automorphism $\alpha$ such that $\phi\circ\alpha$ is a reduced form of $\phi$. By Lemma \ref{l400} and Lemma \ref{l300}, the reduction can be done continuously till some component is reduced to a linear form of $y$. Hence after left and right combine some automorphism, $\phi$ is transferred to be of the form 
$\phi'=(y, xh_3)$. Left combine the automorphism $\beta=(y,x)$, and then $\beta\circ\phi'=(x,yh'_3)$ where $h'_3\not=0$. Again by Lemma \ref{l200} and Lemma \ref{l300}, if $h'_3$ is not a constant, then $x$ appears in $h'_3$ and hence $v(yh'_3)$ and $x$ are algebraically dependent. This is impossible, so $h'_3$ has to be a non-zero constant which implies that $\beta\circ\phi'$ is an automorphism. And then $\phi$ is an automorphism.

\

\noindent {\bf Acknowledgement.}  The author is grateful to Jie-Tai Yu for his helpful discussion and supervision during the period of the work being done.

\end{document}